\documentclass[]{article}
\usepackage[top=1in,left=1in,right=1in,bottom=1.5in]{geometry}
\usepackage{latexsym}
\usepackage{epsfig, ecltree, epic, eepic}
\usepackage{enumerate,amsmath,amssymb,dsfont,pifont,amsthm}
\usepackage[dvips]{color}
\usepackage{graphicx}

\numberwithin{equation}{section}

\theoremstyle{plain}
\newtheorem{Thm}{Theorem}[section]
\newtheorem{Cor}[Thm]{Corollary}
\newtheorem{Lem}[Thm]{Lemma}
\newtheorem{Prop}[Thm]{Proposition}

\theoremstyle{definition}
\newtheorem{Def}[Thm]{Definition}

\theoremstyle{remark}

\newtheorem{Ex}[Thm]{Example}

\renewcommand{\leq}{\leqslant}
\renewcommand{\geq}{\geqslant}

\newcommand{\bbc}{\mathbb{C}}
\newcommand{\bbr}{\mathbb{R}}
\newcommand{\bbq}{\mathbb{Q}}
\newcommand{\bbp}{\mathbb{P}}
\newcommand{\bbe}{\mathbb{E}}

\newcommand{\bbn}{\mathbb{N}}
\newcommand{\bbb}{\mathbb{B}}

\newcommand{\cb}{\mathcal{B}}
\newcommand{\cf}{\mathcal{F}}
\newcommand{\cm}{\mathcal{M}}
\newcommand{\cs}{\mathcal{S}}
\newcommand{\ca}{\mathcal{A}}
\newcommand{\cv}{\mathcal{V}}
\newcommand{\cl}{\mathcal{L}}
\newcommand{\cp}{\mathcal{P}}

\newcommand{\pfe}{\longmapsto}
\newcommand{\pf}{\longrightarrow}
\newcommand{\uparrows}{\uparrow\uparrow}

\newcommand{\abs}[1]{\left| #1 \right|}
\newcommand{\norm}[1]{\left\| #1 \right\|}

\newcommand{\loi}[2]{\left] #1 , #2 \right]}                             
\newcommand{\roi}[2]{\left[ #1 , #2 \right[}

\newcommand{\ci}[2]{\left[ #1, #2 \right]}

\DeclareMathOperator{\Var}{Var}

\begin{document}

\allowdisplaybreaks

\title{\bfseries On the Semimartingale Nature of Feller Processes with Killing}

\author{%
    \textsc{Alexander Schnurr}%
    \thanks{Lehrstuhl IV, Fakult\"at f\"ur Mathematik, Technische Universit\"at Dortmund,
              D-44227 Dortmund, Germany, Tel: +49-231-755-5917, Fax: +49-231-755-3064, e-mail:
              \texttt{alexander.schnurr@math.tu-dortmund.de}
              }
    }

\date{}

\maketitle

\begin{abstract}
Let $U$ be an open set in $\bbr^d$. We show that under a mild assumption on the richness of the generator a Feller process in $(U,\cb(U))$ with (predictable) killing is a semimartingale. To this end we generalize the notion of semimartingales in a natural way to those `with killing'. Furthermore we calculate the semimartingale characteristics of the Feller process explicitly and analyze their connections to the symbol. Finally we derive a probabilistic formula to calculate the symbol of the process.

\vspace*{5mm}


\vspace*{5mm}

\end{abstract}

\section{Introduction}

Our aim is to show that every `rich' Feller process (with predictable killing) is a semimartingale and to further investigate its structural properties. To this end we analyze the connection between the semimartingale characteristics and the symbol of the process. For the latter concept we derive a stochastic formula.

Some interesting but slightly technical criteria when a strong Markov process (on $\bbr^d$, without killing) is a semimartingale are given in \cite{vierleute}. We use a different approach and we prove along the way several estimates which are of interest on their own right. Our starting point is a classical result of the theory of pseudo-differential operators which is due to Courr\`ege \cite{courrege}. Let us emphasize that the case on $\bbr^d$ without killing is always included in our considerations.

For Feller processes on $\bbr^d$ satisfying the so called `growth condition' (formula \eqref{growthcond} below) it is known that several path properties of the process can be derived by analyzing the symbol (cf. \cite{schilling98}, \cite{schillingschnurr}). Therefore the present paper can (and should) be seen as a first step to generalize the respective results to processes with killing, defined on subsets of $\bbr^d$ which do not have to fulfill the `growth condition'. This is part of ongoing research and not in the scope of the present paper.

Let us furthermore emphasize that boundary conditions are not taken into account here. Cases in which the process hits the boundary are a priori excluded by our definitions. The interplay between boundary conditions and the symbol are definitely an interesting topic but it would lead us into a totally different direction and would go beyond the scope of the present paper. 

As there are different conventions in defining Feller processes in the literature, we fix some terminology and several notations in Section 2. While the theory of Feller processes defined on open subsets of $\bbr^d$ with killing is well developed, this is not the case for semimartingales. We establish the theory for this class of processes in Appendix A. Most of the results of the appendix are straight-forward extensions of the classical theory. However, for the readers convenience we decided to include this material and emphasize the differences between these processes and classical semimartingales. Furthermore this appendix preserves us from writing: `By the version of Proposition ... with killing'. In Section 3 we show that every Feller process in the above sense is a semimartingale with killing and even an It\^o process (cf. Definition \ref{def:ito}). Afterwards we introduce the so called stochastic symbol of the process and analyze the relationship between this concept, the generator of the process and the semimartingale characteristics. In order to round out the paper some examples and applications are enclosed in Section 5. The main results are Theorem \ref{thm:fellerinito} and Theorem \ref{thm:symbol}.

Most of the notation we are using is standard: if random variables are considered `=' means a.s. and in the context of stochastic processes it means indistinguishable. We write $\cp$ for the predictable $\sigma$-algebra. For the stochastic integral we use the notations
\[
\int_0^t H_s \, dX_s = (H \bullet X)_t \hspace{10mm} (t\geq 0)
\]
interchangeably. For unexplained notation see \cite{jacodshir}. For the well-known classes  of stochastic processes $\cm, \, \cv, \, \cs  ...$ the killed counterparts are written as $\cm^\dagger, \, \cv^\dagger, \, \cs^\dagger...$ (for details consult the appendix).

\section{Feller Processes and General Framework}

In the context of Feller processes we follow mainly \cite{taira} and \cite{revuzyor}: 
let $U\subseteq \bbr^d$ be an open set and $U_\Delta:=U\cup \{\Delta \}$ its one-point compactification. The jumps of a process with values in $U$ are all contained in the open set
\[
U-U:=\{z\in\bbr^d: \text{there are } x,y\in U \text{ such that } z=y-x \}.
\]
We write $\cb(U),\cb(U_\Delta), \cb(U-U)$ for the respective Borel sets. Vectors in $\bbr^d$ are thought of as column vectors. A transposed vector is denoted by $x'$ and the vector entries by $x^{(1)},...,x^{(d)}$. Furthermore $B_b(U)$ are the bounded, real valued, Borel measurable functions and $C_0(U)$ denotes the continuous, real valued functions which are vanishing at $\Delta$, that is, for every $\varepsilon>0$ there exists a compact set $K\subseteq U$ such that $|u(x)|<\varepsilon$ for every $x\in U \backslash K$. $C_0(U)$ is a closed subspace of the Banach space of continuous bounded functions $(C_b(U),\norm{\cdot})$ and hence itself a Banach space. Here and for the remainder of the present article $\norm{\cdot}$ denotes the supremum norm. Occasionally we will make use of the following conventions: $\Delta \cdot 0=0$, $\Delta + x = \Delta$ and $\Delta - x =\Delta$ for every $x\in U_\Delta$. For a compact set $K\subseteq U$ and $u:U\to \bbr$ we write
\[
\norm{u}_{K}:=\sup_{x\in K} \abs{u(x)}.
\]
As a convention, every real valued function $u$ on $U$ is extended to $\widetilde{u}$ on $U_\Delta$ by setting $\widetilde{u}(\Delta)=0$, if not mentioned otherwise. This allows us to identify $C_0(U)$ with a subspace of $C(U_\Delta)$ as follows
\begin{align} \label{identify}
  C_0(U) = \{ \widetilde{u}\in C(U_\Delta) : \widetilde{u}(\Delta)=0 \}.
\end{align}

The compact unit ball around $0$ is denoted by $\bbb$. A function $\kappa:U-U\to\bbr$ is called cut-off function if it is Borel measurable, with compact support and equal to one in a neighborhood of zero. In this case $h(y):= \kappa(y)\cdot y$ is a truncation function in the sense of \cite{jacodshir}. For technical reasons we fix an $R>0$ such that ${ (2R)\bbb }\subseteq U-U$ and assume $1_{ R\bbb }(y) \leq \kappa(y) \leq 1_{(2R)\bbb}(y)$.


Consider a Markov process in the sense of Blumenthal-Getoor (cf. \cite{blumenthalget})
\[
\textbf{X}:=(\Omega, \cf ,(\cf_t)_{t\geq 0},(X_t)_{t\geq 0},\bbp^x)_{x\in U_\Delta}
\]
with state space $(U_\Delta,\cb(U_\Delta))$, which is time-homogeneous and normal, that is,  $\bbp^x(X_0=x)=1$. Let $P_t$ denote its transition function. As usual we interpret $\Delta$ as the cemetery and have, therefore,
\[
  P_t(\Delta,\{\Delta\}) = 1, P_t(\Delta,U)=0.
\]
We associate a semigroup $(T_t)_{t\geq 0}$ of operators on $B_b(U)$ by setting
\[
  T_t u(x):= \int_{U} u(y) P_t(x,dy)=\bbe^x u(X_t)  \hspace{1cm} (t\geq 0, x\in U).
\]
$T_t$ is for every $t\geq 0$ a contractive, positivity preserving and sub-Markovian operator on this space.
If in addition \newline
\hspace*{1cm}(F1) $T_t:C_0(U) \pf C_0(U)$ for every $t\geq 0$ and \newline
\hspace*{1cm}(F2) $\lim_{t\downarrow 0} \norm{T_tu-u} =0$ for every $u\in C_0(U)$ \newline
we call the semigroup and the associated process $X=(X_t)_{t\geq 0}$ Feller. The stochastic basis $(\Omega,\cf,(\cf_t)_{t\geq 0}, \bbp^x)_{x\in U_\Delta}$ is always in the background.

Analogously one could define a semigroup of Markovian operators on $B_b(U_\Delta)$
\[
   \widetilde{T}_t \widetilde{u}(x):= \int_{U} \widetilde{u}(y) P_t(x,dy)=\bbe^x \widetilde{u}(X_t)  \hspace{1cm} (t\geq 0, x\in U_\Delta).
\]
By \eqref{identify} we have the following equivalence on $C_0(U)$
\[
\widetilde{T}_t(\widetilde{u}) = \widetilde{T_t(u)}
\]
which allows us to switch between the two settings.

The generator $(A,D(A))$ is the closed operator given by
\begin{align} \label{generator}
  Au(x):=\lim_{t \downarrow 0} \frac{T_t u(x) -u(x)}{t} \hspace{1cm}  (u\in D(A))
\end{align}
where $D(A)\subseteq C_0(U_\Delta)$ is the set on which the limit \eqref{generator} exists in strong sense, that is, uniformly in $x\in U$.
The following result is taken from \cite{revuzyor} Proposition VII.1.6.
\begin{Lem} \label{lem:generatormg}
For every $u\in D(A)$ and every $x\in U_\Delta$ the process
\begin{align*}
  M_t^{[u]}&=  u(X_t) -u(x) -\int_0^t A u(X_s) \, ds
\end{align*}
is a $\bbp^x$-martingale.
\end{Lem}
In the sequel we assume that the domain of the generator is sufficiently rich, that is,
\begin{align} \label{rich}
  C_K^\infty(U) \subseteq D(A).
\end{align}
A classical result due to P. Courr\`ege (cf. \cite{courrege}) shows that if \eqref{rich} is fulfilled, $-A$ restricted to the test functions $C_K^\infty(U)$ is a pseudo-differential operator with continuous negative definite symbol. This means $A$ can be written as
\begin{align} \label{pseudo}
  Au(x)= - \int_{\bbr^d} e^{ix'\xi} q(x,\xi) \hat{u}(\xi) \, d\xi \hspace{1cm}  (u\in C_K^\infty(U))
\end{align}
where $\hat{u}=1/(2\pi)^d\int e^{-iy'\xi}u(y) dy$ denotes the Fourier transform and $q:U \times \bbr^d \pf \bbc$ is locally bounded and for fixed $x$ a continuous negative definite function in the sense of Schoenberg (cf. \cite{bergforst}, Chapter II) in the co-variable. In particular it admits a L\'evy-Khintchine representation (with a cut-off function $\kappa$ as above)
\begin{align} \begin{split}\label{lkfxa}
  q(x,\xi)&= a(x)-i \ell(x)'  \xi + \frac{1}{2} \xi'Q(x) \xi \\
      &- \int_{U-U} \left(e^{i y'\xi} -1 - i y'\xi \cdot \kappa(y) \right) \,  N(x,dy)
\end{split}\end{align}
where (for every fixed $x\in U$) $a(x)\in\bbr$, $\ell(x) \in \bbr^d$, $Q(x)$ is a positive semidefinite matrix and $N(x,dy)$ is a transition kernel with the property $N(x,\{0\})=0$ and such that for every $u\in C_K(U)$ the function
\begin{align} \label{Npropone}
  x\mapsto \int_{U-U} \abs{y}^2 u(y) \, N(x,dy) \text{ is measurable and locally bounded.}
\end{align}
Furthermore $N$ is bounded at infinity, that is, for every $x\in U$ there exists a $c(x)>0$ such that
\begin{align}\label{Nproptwo}
N(x,((U-U)\backslash V)-x)\leq c(x) <\infty \text{ for every neighborhood } V \text{ of } x.
\end{align}
Let us remark that the integral in \eqref{lkfxa} is actually defined on $U-x$. Therefore we have in addition
\begin{align}
N(x,(U-U)\backslash (U-x))=0.
\end{align}
We restrict ourselves to the case where $a(x)=0$ in \eqref{lkfxa}. Let us remark that unlike in the case of L\'evy processes this does not mean that there is no killing (cf. examples \ref{ex:levykilling} and \ref{ex:superdrift}). The function $q(x,\xi)$ is called the symbol of the operator $-A$ and as well the symbol of the process. More facts on these relationships can be found in the monograph by N. Jacob \cite{niels1} mostly on the space $\bbr^d$. For some interesting results on the conservativeness of the process see \cite{schilling98pos} (Section 5) and \cite{walterhabil} (Chapter 9). It is a well known fact that there exists a c\`adl\`ag version for every Feller process defined on a stochastic basis which fulfills the usual hypotheses. In the following sections we assume every Feller process to have these properties.

In our investigations we always assume that the process exists. The question for which symbols there exists a corresponding Feller process is a matter of ongoing research. See in this context \cite{hoh00}, \cite{abelskassmann09} and \cite{bjoern} and for a survey of the different construction methods \cite{jacobschilling}.

\section{Feller Semimartingales}

In the sequel $X$ denotes a Feller process with generator $(A,D(A))$ such that $C_K^\infty(U) \subseteq D(A)$ and the restriction of $A$ to the test functions is given by \eqref{pseudo}
where $q:U\times\bbr^d \to \bbc$ is the continuous negative definite symbol with L\'evy  triplet $(\ell(x),Q(x),N(x,dy))$ given by
\begin{align} \begin{split}\label{lkfx}
  q(x,\xi)&= -i \ell(x)'  \xi + \frac{1}{2} \xi'Q(x) \xi
      - \int_{U-U} \left(e^{i y'\xi} -1 - i y'\xi \cdot \kappa(y) \right) \,  N(x,dy).
\end{split}\end{align}
Recall that $1_{ R\bbb } \leq \kappa \leq 1_{(2R)\bbb}$ and $(2R)\bbb\subseteq U-U$.

\begin{Def} A Feller process $X$ is called \emph{Feller process with killing} if it is a process with killing (in the sense of Definition \ref{def:killedprocess}) with respect to every $\bbp^x$ ($x\in U$).
\end{Def}

\emph{Remark:} Although the killing time and the announcing sequence do depend on the starting point $x\in U$ we will write $\zeta$ resp. $\tau_n$ instead of $\zeta^x$ resp. $\tau^x_n$ whenever $x$ is fixed. If not mentioned otherwise, $(\tau_n)_{n\in \bbn}$ always denotes a fixed announcing sequence for $\zeta$.


\begin{Thm} \label{thm:fellersemimg}
Let $X$ be a c\`adl\`ag Feller process on $U$ with killing. Let $(A,D(A))$ be the generator of the process with the property $C_K^\infty(U) \subseteq D(A)$ and $A|_{C_K^\infty(U)}$ as above. Then $X$ is a $d$-dimensional semimartingale with killing with respect to every $\bbp^x$ ($x\in U$).
\end{Thm}


\textbf{Proof:}  In this proof we follow mainly \cite{schilling98}. Let $x\in U$.
There exists a sequence of compact sets $(K_n)_{n\in\bbn}$ and a sequence of open sets $(U_n)_{n\in\bbn}$ such that $K_n \subseteq x+n\bbb$, $K_n\uparrow U$ and
\[
  x\in K_1 \subseteq U_1 \subseteq K_2 \subseteq U_2 ...
\]
Let $j\in \{1,..,d\}$ and $(\phi_n)_{n\in\bbn}\subseteq C_K^\infty(U)$ such that
$\phi_n(x)$ is equal $x^{(j)}$ on $K_n$ and vanishes outside $U_n$. Then, by Lemma \ref{lem:generatormg},
\[
M_t^n:=\phi_n(X_t)-\int_0^t A \phi_n(X_s) \, ds
\]
is a $\bbp^x$-martingale for every $n\in \bbn$. Since $\phi_n\in C_K^\infty(U)$, we have $A\phi_n\in C_0(U)$ which implies that $A\phi_n$ is bounded. It follows that the integral term is of finite variation on compacts. Thus, $\phi_n(X)$ is a classical semimartingale for every $n\in\bbn$.
We have $\phi_n(X_t)=X_t^{(j)}$ as long as $X_t$ is in $K_n$, that is, on the interval $[0,\rho_n[$, where $\rho_n:=\sigma_n \wedge \tau_n$ and $\sigma_n:=\sigma_n^x:=\inf\{t\geq 0: X_t\in U_\Delta \backslash K_n \}$. 
Theorem \ref{thm:prelocal} tells us that $X^{(j)}$ is a semimartingale with killing for $j\in\{1,...,d\}$ and hence $X$ is a $d$-dimensional semimartingale with killing. \hfill $\square$

Now we want to further investigate the semimartingale nature of a Feller process.

In earlier papers on this topic (cf. \cite{schilling98}, \cite{schilling98pos}) the following growth condition is often needed
\begin{align} \label{growthcond}
 \sup_{x\in U} \abs{q(x,\xi)} \leq c \cdot (1 + \abs{\xi}^2) \hspace{10mm} (\xi\in \bbr^d).
\end{align}
We show in the following lemma that a local version of this estimate always holds for the symbols we are dealing with.

\begin{Lem} \label{lem:cndfone}
Let (A,D(A)) be the generator of a Feller process $X$ such that $C_K^\infty(U) \subseteq D(A)$ and
$q(x,\xi)$ be its symbol. In this case the following two (equivalent) conditions are always met:
\begin{enumerate}
\item For every compact set $K\subseteq U$ there is a $c_K\geq 0$ such that
\begin{align} \label{localgrowthcond}
 \sup_{x\in K} \abs{q(x,\xi)} \leq c_K \cdot (1 + \abs{\xi}^2) \hspace{10mm}  (\xi\in \bbr^d).
\end{align}
\item  For every compact set $K\subseteq U$ and $\varepsilon >0$ the three quantities 
$\norm{\ell}_{K},\ \norm{Q}_{K}$ and $\norm{Nb_\varepsilon}_K$ are finite where $Nb_\varepsilon(x)=\int_{U-U} b_\varepsilon (y) \ N(x,dy)$ and
\begin{align} \label{functionb}
b_\varepsilon(y):=\abs{y}^2 1_{\varepsilon \bbb}(y)+1_{ \varepsilon \bbb^c}(y).
\end{align}
\end{enumerate}
 \end{Lem}

\textbf{Proof:}
By Lemma 3.6.21 of \cite{niels1}, the square root of a continuous negative definite function is sub-additive. Thus, $\sqrt{\abs{ q(x,\xi)}} \leq k \sqrt{\abs{ q(x,\xi/k)}}$, and hence,
for every $k\in \bbn$ and $x\in U$
\begin{align*}
\abs{ q(x,\xi)} =k^2 \abs{q\left(x,\frac{\xi}{k}\right)}.
\end{align*}
For a given $\xi$ we take $k_0$ to be the integer in $[\abs{\xi}, \abs{\xi}+1[$. Then,
\[
\abs{ q(x,\xi)}\leq (1+\abs{\xi})^2 \cdot \abs{q\left(x,\frac{\xi}{k_0}\right)}
\leq 2(1+\abs{\xi}^2) \cdot \sup_{\abs{\eta}\leq 1} \abs{q(x,\eta)}.
\]
Let $K\subseteq\bbr^d$ be a compact set. We obtain
\begin{align*}
\sup_{x\in K} \abs{q(x,\xi)}
\leq  \left(2\sup_{x\in K} \sup_{\abs{\eta}\leq 1} \abs{q(x,\eta)} \right) \cdot (1+\abs{\xi}^2) 
\leq c_K \cdot (1+\abs{\xi}^2)
\end{align*}
since the symbol $q(x,\xi)$ is locally bounded (cf. Section 2).

Next we observe that the local growth condition \eqref{localgrowthcond} is equivalent to the local boundedness of the constitutes of the L\'evy triplet $\ell, \, Q$ and $N(\cdot,dy)$. The proof works analogously to \cite{schilling98} Lemma 2.1. In the end one has to use the fact that for every compact set $K\subseteq\bbr^d$ and $\varepsilon>0$
\[
\frac{\abs{y}^2}{1+\abs{y}^2}
\asymp b_\varepsilon(y)
\]
for $y\to 0$ and $\abs{y}\to\infty$ (if $U$ is unbounded). Here $f\asymp g$ means that there exists positive constants $c_1$ and $c_2$ such that 
\[
  c_1\abs{g(y)} \leq \abs{f(y)} \leq c_2 \abs{g(y)}
\]
for $\abs{y}$ small (resp. large) enough.
Thus the result. \hfill $\square$

\emph{Remarks:} (a) In the context of the L\'evy-Khintchine formula the particular truncation function $\abs{y}^2/(1+\abs{y}^2)$ (cf. \cite{sato} Remark 8.4) is often used (see \cite{jacobschillinglkf}). But let us mention that this is not a cut-off function as we have introduced it in Section 2. It turns out that in order to establish a neat representation of the characteristics of a Feller process one should use the same cut-off function for both: the L\'evy-Khintchine representation of the symbol and the semimartingale characteristics.

(b) The second statement of Lemma \ref{lem:cndfone} implies a locally uniform version of \eqref{Nproptwo}.

Studying the domain $D(A)$ of the generator $A$ it is useful to rewrite it in the so called integro-differential-representation. First we need the following estimate:
\begin{Lem}
For a cut-off function $\kappa=\kappa_R$ as described in Section 2, $y\in U$ and $\xi\in\bbr^d$ we have
\begin{align}\label{cndfintesti}
  \abs{e^{iy'\xi} - 1 - iy'\xi \kappa(y)} \leq 2(R+1)(1+\abs{\xi}^2)b_{R\wedge 1}(y).
\end{align}
\end{Lem}

\textbf{Proof:}
Recall that $1_{ R\bbb } \leq \kappa \leq 1_{(2R) \bbb }$. Consider $\abs{\exp(iy'\xi)-1-iy\xi \kappa(y)}$. For $\abs{y}\leq R\wedge 1$, it is dominated by $(1/2) (y'\xi)^2\leq (1/2) \abs{\xi}^2\abs{y}^2$. For $R\wedge 1 \leq \abs{y} \leq 2R$, it is dominated by $2+\abs{\xi}\abs{y} \leq 2+(1+\abs{\xi}^2) \cdot 2R$, and for $\abs{y}>2R$ by 2. Combining these yields the claim. \hfill $\square$

We define the following operator on $C_b^2(U)$:
\begin{align} \begin{split} \label{integrodiff}
I_qu(x)&:=\ell(x) '\nabla u(x) + \frac{1}{2}\sum_{j,k=1}^d \Big(Q^{jk}(x) \partial_j \partial_k u(x)\Big) \\
&+ \int_{U-U} \Big( u(x+y) - u(x) - y'\nabla u(x) \cdot \kappa(y) \Big) N(x,dy).
\end{split}
\end{align}
For $u$ in $C_K^\infty (U)$ it follows from \eqref{pseudo} and \eqref{lkfx} after a change in the order of integration made possible by \eqref{Npropone}, \eqref{Nproptwo} and the estimate \eqref{cndfintesti} that 
\[
A|_{C_K^\infty(U)}= I_q|_{C_K^\infty(U)}
\]
where $\ell,Q$ and $Nb_\varepsilon$ in \eqref{integrodiff} are locally bounded.

In order to get control over the last term of $I_q$ the following estimate is useful.

\begin{Lem}
Let $K\subseteq U$ be a compact set, $u\in C_b^2(U)$ and $\kappa=\kappa_R$ the cut-off function from above.
For $\varepsilon\in ]0,1[$ such that $K+\varepsilon \bbb \subseteq U$, $x\in K$ and $y\in U$ we have
\begin{align} \label{integrandestimateone} \begin{split}
  &\Big|u(x+y)-u(x)-y'\nabla u(x) \cdot \kappa(y)\Big|\\
  &\leq (2R) b_{\varepsilon}(y) \left( \norm{u} + \sum_{\abs{\alpha}=1} \norm{\partial^\alpha u}_{K} + \sum_{\abs{\alpha}=2} \norm{\partial^\alpha u}_{K+\varepsilon \bbb}\right).
\end{split}
\end{align}
In particular it follows
\begin{align} \label{integrandestimatetwo}
  \Big|u(x+y)-u(x)-y'\nabla u(x) \cdot \kappa(y)\Big|
    \leq (2R) b_\varepsilon(y) \sum_{\abs{\alpha}\leq 2} \norm{\partial^\alpha u}.
\end{align}
\end{Lem}

\textbf{Proof:}
Fix $x\in K$ and $y\in\bbr^d$. For $y\in\varepsilon \bbb$ the left-hand side of \eqref{integrandestimateone} is bounded by
\[
\frac{1}{2} \sum_{j,k=1}^d  \norm{\partial_j\partial_k u}_{K+\varepsilon \bbb} \abs{y^{(j)} y^{(k)}} \kappa(y) \leq \frac{1}{2}  \abs{y}^2 \sum_{j,k=1}^d \norm{\partial_j \partial_k u}_{K+\varepsilon \bbb}
\]
where we used a Taylor expansion. For $y$ outside $\varepsilon \bbb$ the same left-hand side is bounded by 
\[
2\norm{u}+\abs{y} \cdot \abs{\nabla u(x)} \cdot \kappa(y) \leq 2 \norm{u} + 2R \abs{\nabla u(x)}
\]
since $\kappa$ vanishes outside $2R\bbb$. Combining these two estimates we obtain the right-hand side of \eqref{integrandestimateone}. Equation \eqref{integrandestimatetwo} follows directly from \eqref{integrandestimateone}.
\hfill $\square$

Now we prove a local version of the inequality (2.9) of \cite{schilling98}. Since we have not demanded that the test functions are a core of $A$ (cf. \cite{ethierkurtz} Section 1.3), there might be different extensions of $A|_{C_K^\infty(U)}$ to a Feller generator, that is, the generator of a Feller process. The results of this section are true for every such extension.

\begin{Lem}
Let $I_q$ be the operator defined in \eqref{integrodiff}. Then for every compact set $K\subseteq U$ there is a constant $d_K>0$ such that for $u\in C_b^2(U)$
\begin{align}\label{insidesupportbetter}
\norm{I_qu}_{K} \leq d_K \cdot \norm{u} + \sum_{\abs{\alpha}\in\{1,2\}} \norm{\partial^\alpha u}_{K+\varepsilon \bbb}
\end{align}
where $\abs{\alpha}=\abs{\alpha^{(1)}}+...+\abs{\alpha^{(d)}}$.
In particular we have
\begin{align} \label{insidesupport}
\norm{I_qu}_{K} \leq d_K \cdot \left(\sum_{\abs{\alpha}\leq 2} \norm{\partial^\alpha u} \right).
\end{align}
\end{Lem}

\textbf{Proof:} Let $K\subseteq U$ be a compact set. We may assume that $\ell=0$ and $Q=0$ as for these `coefficients' the inequality is clear. Uniformly for $x\in K$ we obtain using \eqref{integrandestimateone}
\begin{align*}
\abs{I_qu(x)} \leq 2R \cdot Nb_\varepsilon(x) 
\left(\norm{u} + \sum_{\abs{\alpha}\in\{1,2\}} \norm{\partial^\alpha u}_{K+ \varepsilon \bbb}\right).
\end{align*}
where $Nb_\varepsilon$ is locally bounded by Lemma \ref{lem:cndfone}. This completes the proof.
\hfill $\square$

The next observation helps us to prove that the domain of the Feller generator is quite rich: the generator $A$ maps the test functions (which in our investigations are always in $D(A)$) into $C_0(U)$. This means in particular that for $u\in C_K^\infty(U)$ we have
\[
\lim_{x\to \Delta}Au(x)=0.
\]
For $x$ outside of supp$(u)$ this reads
\begin{align} \label{outofsupport}
\abs{\int_{U-U} u(x+y) \, N(x,dy) } \xrightarrow[x\to \Delta]{}0.
\end{align}

\begin{Thm} \label{thm:richdomain}
Let (A,D(A)) be the generator of a Feller process $X$ such that $C_K^\infty(U) \subseteq D(A)$ Then we have
\[
C_K^2(U) \subseteq D(A).
\]
\end{Thm}

\textbf{Proof:}  The operator $A$ is closed. Let $u\in C_K^2(U)$. Using Friedrichs mollifier we know that there is  a compact set $K\subseteq U$ and a sequence $(u_n)_{n\in\bbn} \subseteq C_K^\infty(U)$ which converges to $u$ in the norm $\sum_{\abs{\alpha}\leq 2} \norm{\partial^\alpha \cdot}$ and has the additional property $\textnormal{supp}(u_n) \subseteq K$ for every  $n\in \bbn$.

If we can show that the sequence $(Au_n)_{n\in\bbn}$ is a Cauchy sequence in $(C_0(U),\norm{\cdot})$ the assertion will follow because of the closedness of the operator. Let $\varepsilon >0$. For $x$ outside of $K$ we have
\[
A(u_m-u_n)(x)=\int_{U-U} (u_m-u_n)(x+y) \, N(x,dy).
\]
Since the convergence of $(u_n)_{n\in\bbn}$ is uniform, the sequence $(u_n)_{n\in\bbn}$ is uniformly bounded (and the support of every $u_n$ is in $K$). Therefore, we can find a non-negative function $f\in C_K^\infty(U)$ such that $-f \leq u_n \leq f$ for every $n\in \bbn$. 
This implies, via \eqref{outofsupport}, that

\[
\int_{U-U} (u_m-u_n)(x+y) \, N(x,dy) \leq 2 \int_{U-U} f(x+y) \, N(x,dy)
\]
and the right-hand side vanishes as $x\to \Delta$. Thus $A(u_m-u_n)$ tends to 0 uniformly in $n,m$, that is, there exists a compact set $\tilde{K} \supseteq K$ such that for $x$ outside $\tilde{K}$ we obtain
\begin{align} \label{outofsupporttwo}
\Big|A(u_m-u_n)(x)\Big| = \abs{\int_{U-U} (u_m-u_n) (x+y) \, N(x,dy)} < \varepsilon.
\end{align}
For $x\in \tilde{K}$ we use formula \eqref{insidesupport}, recalling that $A|_{C_K^\infty(U)}= I_q|_{C_K^\infty(U)}$:
\[
\norm{A(u_m-u_n)}_{\tilde{K}} \leq d_{\tilde{K}} \cdot \sum_{\abs{\alpha}\leq 2}\norm{\partial^\alpha (u_m-u_n)}.
\]
Since $(u_n)_{n\in\bbn}$ converges in the norm $\sum_{\abs{\alpha}\leq 2}\norm{\partial^\alpha \cdot}$ we can find an $N\in \bbn$ such that for every $n,m \geq N$
\[
\norm{A(u_m-u_n)}_{\tilde{K} } < \varepsilon.
\]
Together with formula \eqref{outofsupporttwo} this yields the asserted Cauchy property. \hfill $\square$

Since on $C_K^\infty(U)$ the operators are the same and the image-sequence $(Au_n)_{n\in\bbn}$ in the proof above converges uniformly in $C_0(U)$ we obtain the following.

\begin{Cor} \label{cor:operatoridentity}
Under the assumptions of Theorem \ref{thm:richdomain} we have
\[
A|_{C_K^2(U)} = I_q|_{C_K^2(U)}
\]
where $I_q$ is given by \eqref{integrodiff}.
\end{Cor}

Now we can prove our first main result.

\begin{Thm} \label{thm:fellerinito}
Let $(A,D(A))$ be the generator of a Feller process $X$ such that $C_K^\infty (U) \in D(A)$ and with symbol $q:U\times \bbr^d \to \bbc$ given by \eqref{lkfx}. In this case $X$ is an It\^o process with killing (cf. Definition \ref{def:ito}) and its semimartingale characteristics $(B,C,\nu)$ with respect to the cut-off function $\kappa$ are
\begin{align} \begin{split}
B_t^{(j)}(\omega)&=\int_0^t  \ell^{(j)}(X_s(\omega)) \, dF_s \\
C_t^{jk}(\omega)&= \int_0^t Q^{jk}(X_s(\omega)) \, dF_s \\
\nu(\omega;ds,dy)&=N(X_s(\omega),dy) \, dF_s
\end{split}
\end{align}
for every $\bbp^x, (x\in U)$ where $(\ell, Q, N(\cdot,dy))$ is the triplet which appears in the symbol of the Feller process and $F_s=s\cdot 1_{\roi{0}{\zeta}}(s)+ \Delta \cdot 1_{\roi{\zeta}{\infty}}(s)$.
\end{Thm}

\textbf{Proof:}   We already know that $X$ is a semimartingale with killing. By Theorem \ref{thm:charchar} it suffices to show that for every $u:\bbr^d \to \bbr\cup \{\gamma\}$ such that $u|_{U}\in C_b^2(U)$, $u(\Delta)=\gamma$ and $u|_{U_\Delta^c}=0$ the process given by
\begin{align*}
\widetilde{M}_t^{[u]}&= u(X_t)-u(X_0)-\sum_{j=1}^d \int_0^{t}  \Big(\partial_j u (X_{s-})\ell^{(j)}(X_{s-}) \Big)\, dF_s \\
&- \frac{1}{2} \sum_{j,k=1}^d \int_0^{t}  \Big(\partial_j\partial_k u (X_{s-}) Q^{jk}(X_{s-}) \Big)\, dF_s \\
&- \int_0^t \int_{U-U} \Big(u(X_{s-}+y) -u(X_{s-}) - \kappa(y) y'\nabla u(X_{s-})\Big) \, N(X_{s-},dy) \, dF_s
\end{align*}
is a local martingale with killing for every $\bbp^x \ (x\in U)$. Let $u$ be such a function and $x\in U$. We write $\zeta=\zeta^x$ for the predictable killing time with respect to $\bbp^x$ and $(\tau_n)_{n\in\bbn}$ for its announcing sequence. Now let $(U_m)_{m\in\bbn}$ and $(K_m)_{m\in\bbn}$ be the sequences of sets defined in the proof of Theorem \ref{thm:fellersemimg}. Furthermore let $(\kappa_m)_{m\in\bbn}\subseteq C_K^\infty(\bbr^d_\Delta)$ be a sequence of (smooth) cut-off functions such that
\[
  \kappa_m=1 \text{ on } K_m, \ \kappa_m= 0 \text{ on } U_m^c \text{ and } \kappa_m\uparrow 1_U
\]
The sequence $(u_m)_{m\geq 0}$ defined by $u_m:= u \cdot \kappa_m$ is bounded uniformly in $m$.
By definition we have $u_m \to u$, $\partial_j u_m \to \partial_j u$ and $\partial_j\partial_k u_m \to \partial_j \partial_k u$ for $j,k=1,...,d$ where the convergence is locally uniform on $U$. Furthermore, we have $u_m|_{U_\Delta} \in C_K^2(U_\Delta)\subseteq D(A)$.
By Lemma \ref{lem:generatormg} we conclude that
\begin{align*}
M_t^{[u_m]}&= u_m(X_t)-u_m(X_0)-\int_0^t Au_m(X_s) \, ds \\
&= u_m(X_t)-u_m(X_0)-\int_0^t Au_m(X_{s-}) \, ds
\end{align*}
is a $\bbp^x$-martingale. Using the explicit representation of $I_q$ (see \eqref{integrodiff}) we obtain by Corollary \ref{cor:operatoridentity}
\begin{align*}
M_t^{[u_m]} = \widetilde{M}_t^{[u_m]} \hspace{10mm} \text{ on } \roi{0}{\zeta}.
\end{align*}
Let $\sigma_n:=\sigma_n^x:=\inf\{ t\geq 0 :  X_t \in U_\Delta \backslash K_n \}$ for every $n\in\bbn$ and $\rho_n:= \sigma_n \wedge \tau_n$. Then $(\rho_n)_{n\in\bbn}$ is a sequence of stopping times such that $\rho_n \uparrows \zeta$. The stopped processes $(\widetilde{M}_t^{[u_m]})^{\rho_n}$ are martingales for all $m,n\in\bbn$. Let us have a closer look at these processes:
\begin{align*}
(\widetilde{M}_t^{[u_m]})^{\rho_n}&= u_m(X_t^{\rho_n})-u_m(X_0^{\rho_n})-\int_0^{t\wedge \rho_n} I_q u_m(X_{s-}) \, ds \\
&= u_m(X_t^{\rho_n})-u_m(X_0^{\rho_n})-\int_0^{t} I_q u_m(X_{s-}) 1_{\ci{0}{\rho_n}} (s) \, ds. \\
\end{align*}
Using again the explicit representation of $I_q$ we can write:
\begin{align} \begin{split} \label{mainformula}
\hspace{-10mm}(\widetilde{M}_t^{[u_m]})^{\rho_n}&= u_m(X_t^{\rho_n})-u_m(X_0^{\rho_n}) \\
&-\sum_{j=1}^d \int_0^{t}  \Big(\partial_j u_m (X_{s-})\ell^{(j)}(X_{s-}) 1_{\ci{0}{\rho_n}}(s) \Big)\, ds \\
& - \frac{1}{2} \sum_{j,k=1}^d \int_0^{t}  \Big(\partial_j\partial_k u_m (X_{s-}) Q^{jk}(X_{s-}) 1_{\ci{0}{\rho_n}}(s) \Big)\, ds \\
&- \int_0^t \int_{U-U} \Big((u_m(X_{s-}+y) -u_m(X_{s-}) - \kappa(y) y'\nabla u_m(X_{s-}))\\
&\hspace{20mm} \times 1_{\ci{0}{\rho_n}}(s) \Big)\, N(X_{s-},dy) \, ds.
\end{split}
\end{align}

Since $(\widetilde{M}_t^{[u_m]})^{\rho_n}$ is a martingale for every $m,n\in\bbn$ we obtain for $r\leq t$ and  $F\in \cf_r$
\begin{align*}
\int_F (\widetilde{M}_t^{[u_m]})^{\rho_n} \, d\bbp^x = \int_F (\widetilde{M}_r^{[u_m]})^{\rho_n} \, d\bbp^x.
\end{align*}
If we show that for every $r\leq t$ and  $F\in \cf_r$ (the case $r=t$ is included)
\begin{align} \label{laststep}
\int_F (\widetilde{M}_t^{[u_m]})^{\rho_n} \, d\bbp^x \xrightarrow[m\to\infty]{} \int_F (\widetilde{M}_t^{[u]})^{\rho_n} \, d\bbp^x
\end{align}
we will obtain that $(\widetilde{M}_t^{[u]})^{\rho_n}_{t\geq 0}$ is a martingale which will yield in turn that $(\widetilde{M}_t^{[u]})_{t\geq 0}$ is a local martingale with killing and hence the result.
Therefore, the only thing which remains to be proved is \eqref{laststep}. \newline We fix $n\in \bbn$, $r\leq t$ and  $F\in \cf_r$ and show the convergence separately for every term in \eqref{mainformula}.
We start with term number one: the sequence $(u_m)_{m\in\bbn}$ is uniformly bounded. Furthermore, we have that $u_m(X_t^{\rho_n})$ converges pointwise to $u(X_t^{\rho_n})$ and by dominated convergence we obtain
\[
\int_F u_m(X_t^{\rho_n}) \, d\bbp^x \xrightarrow[m\to\infty]{} \int_F u(X_t^{\rho_n}) \, d\bbp^x.
\]
The second term, $u_m(X_0^{\rho_n})$, works alike. \newline
In term three we obtain (for $j=1,...,d$) that $\ell^{(j)}(X_{s-}) \cdot 1_{\ci{0}{\rho_n}}(s)$ is bounded because $\ell$ is locally bounded (see Lemma \ref{lem:cndfone}) and $X_{s-}\in K_n$ on $\ci{0}{\rho_n}$. Furthermore, $\partial_j u_m(X_{s-})$ coincides with $\partial_j u(X_{s-})$ for $m\geq n+1$. Using the dominated convergence theorem on the space $F\times [0,t]$ we obtain
\begin{align*}
&\int_F\int_0^t  \partial_j u_m(X_{s-}) \ell^{(j)}(X_{s-}) 1_{\ci{0}{\rho_n}}(s) \, ds \, d\bbp^x \\ &\xrightarrow[m\to\infty]{} \int_F\int_0^t  \partial_j u(X_{s-}) \ell^{(j)}(X_{s-}) 1_{\ci{0}{\rho_n}}(s) \, ds \, d\bbp^x.
\end{align*}
The fourth term works like the third one. The only difference is that now second derivatives are used.\newline
The integrand in term five
\[
u_m(X_{s-}(\omega)+y) -u_m(X_{s-}(\omega)) - \kappa(y) y'\nabla u_m (X_{s-}(\omega))\cdot 1_{\ci{0}{\rho_n(\omega)}}(s)
\]
coincides for $m\geq n+1$ with
\[
u(X_{s-}(\omega)+y) -u(X_{s-}(\omega)) - \kappa(y) y'\nabla u (X_{s-}(\omega))\cdot 1_{\ci{0}{\rho_n}(\omega)}(s)
\]
The integrability is this time given by \eqref{Npropone} and by \eqref{integrandestimateone} where $K=K_n$ and $\varepsilon$ is chosen such that $K+ \varepsilon \bbb\subseteq K_{n+1}$.
We have thus established \eqref{laststep} and the result follows.
\hfill $\square$

\section{The Symbol of a Feller Semimartingale}

A stochastic formula for the symbol was first analyzed by N. Jacob in \cite{nielsursprung}. R. L. Schilling has proved a more general result in the context of Feller processes (\cite{schilling98pos}). Only the growth condition \eqref{growthcond} is needed there. In the recent paper \cite{schillingschnurr} we have shown that even the growth condition is not necessary. Here we generalize the results to the present setting of a Feller semimartingale with killing.

We write for $\xi\in\bbr^d$
\[
e_\xi(x):=\begin{cases}e^{i x'\xi}& \text{if } x\in U \\ 0 &\text{if } x=\Delta. \end{cases}
\]

\begin{Def}
    Let $X$ be a $U$-valued Markov process with killing which is normal. Fix a starting point $x$ and let $K\subseteq U$ be a compact neighborhood of $x$. Define $\sigma$ to be the first exit time of $X$ from $K$:
    \begin{gather} \label{stopping}
        \sigma:=\sigma^x_K:=\inf\big\{t\geq 0 : X_t\in U_\Delta \backslash K  \big\}.
    \end{gather}
    The function $p:U\times \bbr^d \rightarrow \bbc$ given by
    \begin{gather} \label{symbol}
         p(x,\xi):= -\lim_{t\downarrow 0}  \frac{\bbe^x\Big(e_\xi(X_t^\sigma-x)  - 1\Big)}{t}   = -\lim_{t\downarrow 0}  \frac{\int_{\{X_t^\sigma \neq \Delta \}}e_\xi(X_t^\sigma-x) \, d\bbp^x- 1}{t}
    \end{gather}
    is called the \emph{(probabilistic) symbol of the process}, if the limit exists for every $x\in U$, $\xi\in\bbr^d$ independently of the choice of $K$.
\end{Def}

\emph{Remarks:}
(a) For fixed $x$ the function $p(x,\xi)$ is negative definite as a function of $\xi$. This can be shown as follows: for every $t > 0$ the function $\xi\mapsto \bbe^x e_\xi(X_t^\sigma-x)$ is the characteristic function of the sub-probability measure $\mu(X_t^\sigma-x\in dy)$ where $\mu:=\bbp^x|_{\{X_t^\sigma\neq \Delta\}}$. Therefore it is a continuous positive definite function. By Corollary 7.7 of \cite{bergforst} we conclude that $\xi\mapsto -\bbe^x(e_\xi(X_t^\sigma-x)  - \mu(\Omega))$ is a continuous negative definite function. Since the negative definite functions are a cone which contains the positive constants and which is closed under pointwise limits (cf. \cite{niels1} Lemma 3.6.7), the assertion is proved. Note, however, that $\xi\mapsto p(x,\xi)$ is not necessarily continuous.

(b) In \cite{mydiss} the following is shown for the class of It\^o processes on $\bbr^d$ without killing: fix $x\in\bbr^d$; if the limit \eqref{symbol} exists for one compact neighborhood of $x$, then it exists for every compact neighborhood of $x$ and the respective limits coincide.

If $X$ is a Feller process satisfying \eqref{rich} the symbol $p(x,\xi)$ is exactly the negative definite symbol $q(x,\xi)$ which appears in the pseudo-differential representation of its generator \eqref{pseudo}. A posteriori this justifies the name. In order to obtain this we need the following result which is an immediate consequence of Lemma \ref{lem:generatormg}.
\begin{Lem} \label{lem:dynkin}
Let $X$ be a Feller process with generator $A$. Let $\sigma$ be a stopping time. Then we have
\begin{gather} \label{dynkin}
    \bbe^x \int_0^{\sigma \wedge t} A u(X_s) \,ds = \bbe^x u(X_{\sigma \wedge t})-u(x)
\end{gather}
for all $t>0$ and $u\in D(A)\subseteq C_0(U_\Delta)$.
\end{Lem}

Our second main result gives a formula to directly calculate the symbol, without even writing down the semigroup or the generator:

\begin{Thm}\label{thm:symbol}
Let $(A,D(A))$ be the generator of a Feller process $X$ such that $C_K^\infty (U) \subseteq D(A)$ and with symbol $q:U\times \bbr^d \to \bbc$ given by \eqref{lkfx}.
Let $\sigma$ be as defined in \eqref{stopping}. 
If $x\mapsto q(x,\xi)$ is finely continuous (see \cite{blumenthalget} Section II.4) for each $\xi\in\bbr^d$, then the probabilistic symbol $p$ exists and coincides with the symbol $q$ of the generator.
\end{Thm}

\emph{Remarks:} (a) Combining the result above with Theorem \ref{thm:fellerinito} provides a nice approach to directly calculate the semimartingale characteristics of a Feller process.

(b) Having calculated the symbol, one can write down the generator of the process using formula \eqref{pseudo}.

(c) Let us emphasize that fine continuity is a weaker condition than ordinary continuity.
Even the assumption that $x\mapsto q(x,\xi)$ is continuous would not be a severe restriction. All non-pathological known examples of Feller processes satisfy this condition.

\textbf{Proof:}
Fix $x\in U$ and $\xi\in\bbr^d$.
Let $K_1:=K$ and $(U_n)_{n\in\bbn}$ and $(K_n)_{n\in\bbn}$ be sequences of sets such that $K_n\uparrow U$ and
\[
  x\in K_1 \subseteq U_1 \subseteq K_2 \subseteq U_2 ...
\]
where the $K_n$ are compact and the $U_n$ are open sets for every $n\in\bbn$. Furthermore let $(\kappa_n)_{n\in\bbn}\subseteq C_K^\infty(U)$ be a sequence of (smooth) cut-off functions such that
\[
  \kappa_n=1 \text{ on } K_n, \ \kappa_n= 0 \text{ on } U_n^c \text{ and } \kappa_n\uparrow 1_U
\]
The sequence $(u_n)_{n\geq 0}$ defined by $u_n:= e_\xi \cdot \kappa_n$ is uniformly bounded in $n$.
By definition we have $u_n = e_\xi$, $\partial_j u_n = \partial_j e_\xi$ and $\partial_j\partial_k u_n = \partial_j \partial_k e_\xi$ for $j,k=1,...,d$ on $K_2$ for $n\geq 3$. Furthermore, we have $u_n \in C_K^2(U)\subseteq D(A)$.

By the bounded convergence theorem and Dynkin's formula \eqref{dynkin} we see
\begin{align*}
    \bbe^x \Big(e_\xi(X_t^\sigma-x)  - 1\Big)
    &=\lim_{n\to\infty} \big(\bbe^x  u_n(X_t^\sigma) e_{-\xi}(x) -1 \big) \\
    &=\lim_{n\to\infty} \big(\bbe^x \kappa_n (X_t^\sigma) e_\xi(X_t^\sigma) e_{-\xi}(x) -1 \big) \\
    &=e_{-\xi}(x) \lim_{n\to\infty} \bbe^x \big(\kappa_n(X_t^\sigma) e_\xi (X_t^\sigma) - \kappa_n(x) e_\xi(x) \big)\\
    &=e_{-\xi}(x) \lim_{n\to\infty} \bbe^x \int_0^{\sigma \wedge t} A(\kappa_n  e_\xi)(X_s) \,ds \\
    &=e_{-\xi}(x) \lim_{n\to\infty} \bbe^x \int_0^{\sigma \wedge t} A(\kappa_n  e_\xi)(X_{s-}) \,ds.
\end{align*}
The last equality follows since we are integrating with respect to Lebesgue measure and since a c\`adl\`ag process has a.s.\ a countable number of jumps. By the particular choice of the sequence $(\kappa_n)_{n\in\bbn}$ and by \eqref{insidesupportbetter}, where one has to take $\varepsilon < dist(K_1,U_1^c)$, we can use again the dominated convergence theorem to pass the limit inside the integral and note that $\lim_{n\to\infty} A(\kappa_n e_\xi) = e_\xi p(\cdot,\xi)$ by a classical results due to P. Courr\`ege (cf. \cite{courrege} Sections 3.3 and 3.4). Thus
\begin{align*}
    \lim_{t\downarrow 0} \frac{\bbe^x \Big(e_\xi(X_t^\sigma-x)  - 1\Big)}{t}
    &=-e_{-\xi}(x) \lim_{t \downarrow 0} \bbe^x\left( \frac 1t \int_0^t e_\xi(X_{s-}) p(X_{s-},\xi) 1_{\roi{0}{\sigma}}(s) \,ds\right)
\end{align*}
and we may replace $X_{s-}$ by $X_s$ because we are integrating with respect to Lebesgue measure. The process $X$ is bounded on $\roi{0}{\sigma}$. On the same stochastic interval the mapping $s\mapsto p(X,\xi)$ is bounded and it is right continuous for every $\xi\in\bbr^d$ by the fine continuity of $x\mapsto q(x,\xi)$ (cf. \cite{blumenthalget} Theorem 4.8). Since we have
\[
 \frac{1}{t} \int_0^t e_\xi(X_{s}) p(X_{s},\xi) 1_{\roi{0}{\sigma}}(s) \,ds
 \leq \sup_{0\leq s\leq t}  \Big|e_\xi(X_{s}) p(X_{s},\xi) 1_{\roi{0}{\sigma}}(s)-e_\xi(x)q(x,\xi)\Big|
\]
we obtain by the bounded convergence theorem
\begin{gather*}
    \lim_{t\downarrow 0} \frac{\bbe^x \Big(e_\xi(X_t^\sigma-x)  - 1\Big)}{t}
    =-e_{-\xi}(x) e_{\xi}(x) q(x,\xi)
    =-q(x,\xi).
\end{gather*}
\hfill $\square$

\section{Some Examples and Applications}

\begin{Ex} \label{ex:levy}
Let $X$ be a $d$-dimensional \emph{L\'evy process}. In this case the characteristic function can be written in the following way:
\[
  \bbe^x \left( e^{i \xi'(X_t-x)} \right) = \bbe^0 \left( e^{i \xi'X_t} \right) = e^{-t\cdot \psi(\xi)}
\]
where $\psi:\bbr^d\to\bbc$ is a continuous negative definite function in the sense of Schoenberg. In this case we have for every $x\in\bbr$
\[
\psi(\xi)=p(x,\xi)=q(x,\xi)
\]
The same result holds true for a \emph{subordinator} defined on $]0,\infty[$.
\end{Ex}

\begin{Ex} \label{ex:levykilling}
It is well known that a symbol of the form
\begin{align*}
q(x,\xi)= a -i \ell'  \xi + \frac{1}{2} \xi'Q \xi
      + \int_{U-U} \left(e^{i y'\xi} -1 - i y'\xi \cdot \kappa(y) \right) \,  N(dy)
\end{align*}
gives rise to a \emph{L\'evy process with killing}, that is, plugging this symbol into formula \eqref{pseudo} one obtains an operator which has an extension to the generator of such a process (cf. \cite{nielsold} Chapter 2). The killing time $\zeta$ of such a process is not predictable. Therefore, it does not fit into the setting of Section 3. However one can still use formula \eqref{symbol} in order to calculate the symbol.
\end{Ex}

\begin{Ex} \label{ex:superdrift}
The \emph{superdrift} (starting in x) is the deterministic Markov process given by
\[
X_t^x=\begin{cases} \frac{1}{\frac{1}{x} -t}& \text{if } t\in [0,1/x[ \\ \Delta &\text{else} \end{cases}
\]
for $x>0$. It is easy to see that this process is an It\^o process with killing and even a Feller process. The symbol of this process is $p(x,\xi)= -ix^2\xi$ and therefore its first characteristic is $B_t=\int_0^t X_s^2 \, ds$.
\end{Ex}

\begin{Ex}\label{ex:spacedepdrift}
Let $X_t=X_0+\int_0^t f(X_s) \ ds$ where $f(x):=\text{sgn}(x)$. This is an It\^o process. Its transition semigroup is given by setting $T_t u(x)$ equal to $u(x+t)$ if $x>0$, to $u(0)$ if $x=0$, and to $u(x-t)$ if $x<0$. Taking $t=1$, we see that $T_t u(0+)=u(1)$ whereas $T_t u(0) =u(0)$ (cf. \cite{mydiss} Example B.6). Thus, if $u\in C_0(\bbr)$ is chosen such that $u(1)\neq u(0)$, then $t\mapsto T_tu$ fails to be continuous, so, $X$ is not Feller. 
The semimartingale characteristics of this \emph{space dependent drift} are $(B,C,\nu)=(X-x, 0,0)$ and the symbol of the process is
\[
p(x,\xi)= -i \cdot \text{sgn}(x) \xi.
\]
This symbol is not continuous in $x$, but it is finely continuous, since $t\pfe\ell(X_t)$ is right continuous for every $\bbp^x \ (x\in\bbr)$.
\end{Ex}

\begin{Ex} \label{ex:levydrivensde}
Let $(Z_t)_{t\geq 0}$ be an $\bbr^n$-valued L\'evy process with symbol $\psi(\xi)$ and consider the stochastic differential equation
\begin{align*}
  dX_t^x&=\Phi(X_{t-}^x) \,dZ_t \\
  X_0^x&=x,\quad x\in\bbr^d,
\end{align*}
where $\Phi: \bbr^d \to \bbr^{d \times n}$ is Lipschitz continuous and bounded. In \cite{schillingschnurr} it is shown that the solution $X^x$ is a Feller process with symbol
\[
p(x,\xi)=q(x,\xi)=\psi(\Phi(x)'\xi).
\]
If $\Phi$ is not bounded the process is an It\^o process and the formula for the symbol remains valid.
\end{Ex}

\appendix

\section{Appendix: Killed Semimartingales and Their Characteristics}

In the whole section $\zeta:\Omega\to [0,\infty]$ denotes a predictable time (cf. \cite{jacodshir}), the so called  \emph{killing time}. We write $\tau_n \uparrow \zeta$ (for $n\to\infty$ a.s.) if the sequence of stopping times $(\tau_n)_{n\in\bbn}$ is monotone increasing and converges a.s. to $\zeta$. If in addition $\tau_n < \zeta$ (on $\zeta \neq 0$) we write $\tau_n \uparrows \zeta$ (for $n\to\infty$ a.s.). In the latter case $(\tau_n)_{n\in\bbn}$ is called an \emph{announcing sequence} for $\zeta$.

The reader might wonder why only predictable killing is taken into account. This is owed to the fact that we want to use semimartingale techniques in order to work with Feller processes. One of the most important tools in doing this is by using localization. And having an approximating sequence of stopping times makes $\zeta$ automatically predictable (see e.g. \cite{jacodshir} Theorem I.2.14). When local martingales with killing are treated in the literature, the same restriction is put on $\zeta$ (cf. \cite{getoorsharpe72}).

\begin{Def} \label{def:killedprocess}
Let $U$ be an open subset of $\bbr^d$. Let $X$ be a stochastic process on the stochastic basis $(\Omega, \cf ,(\cf_t)_{t\geq 0},\bbp)$ with values in $(U_\Delta,\cb(U_\Delta))$ and let $\zeta$ be a predictable time. $X$ is called a \emph{process with killing} at $\zeta$ if $X\cdot 1_{\roi{0}{\zeta}}$ takes only values in $U$ and $X\cdot 1_{\roi{\zeta}{\infty}}$ is idenitically equal to $\Delta$.
\end{Def}

For the remainder of the section we assume all processes to be adapted to the filtration $(\cf_t)_{t\geq 0}$.

\begin{Def}
Let $V$ be an open subset of $\bbr$. And let $A$, $M$ be processes on $V_\Delta$ with killing.
\begin{enumerate}
  \item $A$ is called \emph{of finite variation} if the paths $s\mapsto X_s(\omega)$ are a.s. of finite variation on every compact interval $[0,t]\subseteq[0,\zeta(\omega)[$.
  \item If $A$ is adapted, c\`adl\`ag, of finite variation and $A_0=0$ then we write $A\in \cv^\dagger$.
  \item If $A$ is adapted, c\`adl\`ag, increasing and $A_0=0$ then we write $A\in \cv^{+,\dagger}$.
  \item We write $\Var(A)$ for the \emph{variation process} of $A$, that is, the process such that $\Var(A)_t(\omega)$ is the total variation of the function $s\mapsto A_s(\omega)$ on the interval $[0,t]\subseteq[0,\zeta(\omega)[$ and $\Var(A)=\Delta$ on $\roi{\zeta}{\infty}$.
  \item $M$ is called a \emph{local martingale} if there exists an announcing sequence $\tau_n \uparrows \zeta$ such that $M^{\tau_n}$ is a c\`adl\`ag local martingale (in the classical sense). We write $M\in\cm_{loc}^\dagger$.
  \item If $M$ is adapted, c\`adl\`ag, a local martingale and $M_0=0$ then we write $M\in \cl^\dagger$
\end{enumerate}
\end{Def}

\emph{Remarks:}
(a) $\cl^\dagger$, $\cv^\dagger$ are vector spaces. \newline
(b) $\cl^\dagger$, $\cv^\dagger$, $\cv^{+,\dagger}$ are stable under stopping
(cf. \cite{jacodshir} Lemma I.1.35). If e.g. $M\in\cl^\dagger$ and $\tau$ is a stopping time then $M^\tau$ is again in $\cl^\dagger$ but with the new killing time $\widetilde{\zeta}=\zeta\cdot 1_{\{\tau \geq \zeta \}} + \infty \cdot 1_{\{\tau < \zeta \}}$. The stopping time $\widetilde{\zeta}$ is again a predictable time by Propositions I.2.10 and I.2.11 of \cite{jacodshir}.

\begin{Def}
\begin{enumerate}
  \item We write $A\in \ca_{loc}^{+,\dagger}$ if $A\in\cv^{+,\dagger}$ and there exists an announcing sequence $\tau_n\uparrows \zeta$ such that $\bbe A_{\tau_n}\leq \infty$ for every $n\in\bbn$.
  \item We write $A\in \ca_{loc}^\dagger$ if $A\in\cv^\dagger$ and there exists an announcing sequence $\tau_n\uparrows \zeta$ such that $\bbe(\Var(A))_{\tau_n}\leq \infty$ for every $n\in\bbn$.
\end{enumerate}
\end{Def}

\begin{Lem}
Let $A$ be a predictable process in $\cv^\dagger$. Then there exists an announcing sequence $\tau_n\uparrows\zeta$ such that $\Var(A)_{\tau_n}\leq n$ a.s. In particular $A\in\ca_{loc}^\dagger$.
\end{Lem}

\textbf{Proof:}
Let $B:=\Var(A)$. Obviously this process is predictable. Define
\[
  \sigma_n:=\inf\{t\geq 0: B_t=\Delta \text{ or } B_t\geq n\}
\]
for every $n\in\bbn$. These are predictable times since $B$ is predictable. Now let $\tau(n,k) \uparrows \sigma_n$ (for $k\to\infty$) be their respective announcing sequences. For every $n\in\bbn$ there exists a $k_n\in\bbn$ such that
\[
  \bbp\Big(\tau(n,k_n) < \sigma_n-\frac{1}{n}\Big) \leq 2^{-n}.
\]
We set $\tau_n:=\sup_{m\leq n} \tau(m,k_m)$. Then $\tau_n < \sup_{m\leq n} \sigma_m = \sigma_n$ a.s. Therefore we obtain $B_{\tau_n} \leq n$ a.s. Furthermore we have $\sigma_n \uparrow \zeta$. By the definition of $\tau(n,k_n)$ and $\tau_n$ we conclude that $\tau_n \uparrows \zeta$.
\hfill $\square$

This proof illustrates that one does not have to change too much by introducing a killing for semimartingales. In the remainder of the section we will not give detailed proofs but we will emphasize the differences to the proofs of the classical theory.

\begin{Prop}
Let $A\in\ca_{loc}^{+,\dagger}$. There exists a predictable process $A^p\in\ca_{loc}^{+,\dagger}$ called the predictable projection (or compensator) of $A$ such that $A-A^p\in \cm_{loc}^\dagger$.
\end{Prop}

\textbf{Proof:}
Use a localization argument and the Doob-Meyer decomposition as in the original proof (cf. \cite{jacodshir} Theorem I.3.17.). \hfill $\square$

\begin{Def}
\begin{enumerate}
\item Let $V$ be an open subset of $\bbr$. Let $X$ be a process in $V_\Delta$ with killing. We call $X$ a \emph{semimartingale with killing} (at $\zeta$) if there exists a decomposition
\[
  X_t=X_0+B_t+M_t \hspace{1cm} (t\geq 0)
\]
such that $(B_t)_{t\geq 0} \in \cv^\dagger$ and $(M_t)_{t\geq 0} \in \cl^\dagger$, both killed at $\zeta$. We write $X\in\cs^\dagger$ for this class.
\item Let $U$ be an open subset of $\bbr^d$. A process $X$ in $U_\Delta$ with killing is called a ($d$-dimensional) semimartingale with killing, if each of its components is a semimartingale with killing. We write $X\in\vec{\cs}^\dagger$.
\end{enumerate}
\end{Def}

\emph{Remarks:} (a) Recall that by convention $\Delta +  \Delta = \Delta$.

(b) We have $\cs\subseteq\cs^\dagger$. Just set $\zeta:=\infty$ and $\tau_n:=n$.

The following theorem is very important for our considerations. Therefore we give a proof.

\begin{Thm} \label{thm:prelocal}
Let $X$ be a process with killing at $\zeta$. $X$ is a semimartingale with killing iff there exists a sequence of stopping times $(\tau_n)_{n\in\bbn}$ with $\tau_n\uparrows \zeta$ and a sequence $(Y(n))_{n\in\bbn}$ of semimartingales (in the classical sense)  such that $X=Y(n)$ on $\roi{0}{\tau_n}$.
\end{Thm}

\textbf{Proof:}
Sufficiency is clear. Let $(\tau_n)_{n\in\bbn}$ and $(Y(n))_{n\in\bbn}$ as in the theorem. Since $\tau_n\uparrows \zeta$ we obtain that $X$ is c\`adl\`ag. The process
\[
  Z(n):= (Y(n))^{\tau_n} + \Big(X_{\tau_n}-Y_{\tau_n}(n)\Big) 1_{\roi{\tau_n}{\infty}}
\]
is a semimartingale and admits therefore a decomposition (for every $n\in\bbn$):
\[
  Z(n)=X_0+M(n)+A(n)
\]
where $M(n)\in\cl$ and $A(n)\in\cv$. We have by definition $X^{\tau_n}=Z(n)$. Therefore by setting $M:=\sum_{n\in\bbn} M(n) 1_{\loi{\tau_{n-1}}{\tau_n}}$ (with the convention $\tau_0=0$) and $A:=X-X_0-M$ we obtain on $\roi{0}{\zeta}$
\[
  X=X_0+M+A.
\]
$M\in\cl^\dagger$ since
\[
  M^{\tau_n}=\sum_{1\leq k \leq n} \Big( (M(k))^{\tau_k}-(M(k))^{\tau_{k-1}}\Big)\in\cl
\]
and $A\in \cv^\dagger$ since $A=\sum_{n\in\bbn} A(n) 1_{\loi{\tau_{n-1}}{\tau_n}}$. Hence $X\in\cs^\dagger$.
\hfill $\square$

\begin{Def}
A (one-dimensional) semimartingale $X$ with killing is called a \emph{special semimartingale with killing} if it admits a decomposition
\begin{align} \label{specialdecomp}
X_t=X_0+A_t+M_t\hspace{1cm} (t\geq 0)
\end{align}
such that $M\in\cl^\dagger$ and $A$ is a predictable process in $\cv^\dagger$. We write $\cs_p^\dagger$ for this class of processes.
\end{Def}

\begin{Prop} \label{prop:charspecial}
Let $X\in\cs^\dagger$. The following are equivalent:
\begin{itemize}
  \item[(i)] $X\in\cs_p^\dagger$,
  \item[(ii)] There is a decomposition $X=X_0+A+M$ as in \eqref{specialdecomp} with $A\in\ca_{loc}^\dagger$.
  \item[(iii)] For every decomposition $X=X_0+A+M$ as in \eqref{specialdecomp} we have $A\in\ca_{loc}^\dagger$.
  \item[(iv)] We have $Y_t=\sup_{s\leq t} \abs{X_s-X_0} \in \ca_{loc}^{+,\dagger}$.
\end{itemize}
\end{Prop}

\textbf{Proof:} Cf. \cite{jacodshir} Proposition I.4.23 \hfill $\square$.

\begin{Lem} \label{lem:boundedjumps}
If $X\in\cs^\dagger$ and there exists an $a>0$ such that $\abs{\Delta X}\leq a$ on $\roi{0}{\zeta}$ then $X$ is a special semimartingale with killing.
\end{Lem}

\emph{Remark:} Here and in the sequel we use the common notation $\Delta X:=X-X_{-}$. It should be clear from the context if the cemetery state or the jumps of the process are meant.

\textbf{Proof:}
Let $\sigma_n:= \inf \{ t\geq 0 :  X_t= \Delta \text{ or } \abs{X_t-X_0} >n \}$ and $\rho_n:=\sigma_n \wedge \tau_n$. Then the sequence of stopping times $(\rho_n)_{n\in\bbn}$ converges to $\zeta$ and
\[
\sup_{s\leq \rho_n} \abs{X_s-X_0} \leq n+a
\]
Thus (iv) of Proposition \ref{prop:charspecial} is met.
\hfill $\square$

\emph{Remark:} This shows that if $U$ is bounded we have $\cs^\dagger = \cs_p^\dagger$.

The stochastic integral with killing is defined as follows: let $X$ be a semimartingale with killing at $\zeta^X$ and $H$ be a predictable process with killing at $\zeta^H$. We write $\zeta:=\zeta^H\wedge \zeta^X$ and denote by $(\tau_n)_{n\in\bbn}$ the announcing sequence of $\zeta$ given by $\tau_n:=\tau_n^H\wedge\tau_n^X$. We set for every $n\in\bbn$ on $\roi{0}{\zeta}$
\[
  (H\bullet X)^{\tau_n} := H\bullet X^{\tau_n}
\]
and $(H\bullet X):=\Delta$ on $\roi{\zeta}{\infty}$. 

The process $H\bullet X$ is well defined since for $m\geq n$
\begin{align*}
(H\bullet X^{\tau_m})^{\tau_n}
&= (H 1_{\ci{0}{\tau_n}}) \bullet X^{\tau_m} \\
&= H \bullet(1_{\ci{0}{\tau_n}} \bullet X^{\tau_m}) \\
&= H \bullet X^{\tau_n}.
\end{align*}
By a similar argument it is easy to see that this definition does not depend on the choice of the announcing sequence $(\tau_n)_{n\in\bbn}$. See in this context \cite{maisonneuve} (Section 5) where only the case of continuous semimartingales is treated.

In the Section 3 we need It\^o's formula for processes with killing. In order to obtain this we use the classical It\^o's formula on the stopped processes $X^{\tau_n}$.

\begin{Thm}
Let $\bbr\cup \{ \gamma \}$ be the one-point compactification of $\bbr$. Let $f:U\cup\{\Delta\}\to \bbr \cup \{\gamma\}$ which is twice continuously differentiable on $U$ and such that $f(U)\subseteq \bbr$ and $f(\Delta)=\gamma$. Let $X\in\cs^\dagger$ (with $\zeta$ and $\tau_n$). Then $f(X)$ is a semimartingale with killing and we have for every $n\in\bbn$
\begin{equation} \label{itoformula}
\begin{split}
f(X^{\tau_n})-f(X_0)&=\sum_{j\leq d} \partial_j f(X_-^{\tau_n}) \bullet X^{\tau_n,(j)} \\
&+\sum_{s\leq \cdot} \sum_{j\leq d} \partial_j f(X_{s-}^{\tau_n}) \Big(\Delta X_s^{\tau_n,(j)}-\kappa(\Delta X_s^{\tau_n}) (\Delta X_s^{\tau_n})^{(j)}\Big) \\
&+\frac{1}{2} \sum_{k,j\leq d} \partial_j \partial_k f(X_-^{\tau_n}) \bullet \Big< X^{\tau_n, (j), c}, X^{\tau_n, (k), c} \Big>\\
&+\sum_{s\leq \cdot} \Bigg( f(X_s^{\tau_n})-f(X_{s-}^{\tau_n})-\sum_{j\leq d} \partial_j f(X_{s-}^{\tau_n}) \Delta X_s^{\tau_n,(j)} \Bigg)
\end{split}
\end{equation}
and on $\roi{\zeta}{\infty}$ we have $f(X)=\gamma$.
\end{Thm}

Now we introduce the characteristics of a killed semimartingale. Recall that $\kappa$ is a fixed cut-off function. First we set
\begin{align}
\check{X}(\kappa)&:= \sum_{0\leq s \leq \cdot} \Big(\Delta X_s - \kappa(\Delta X_s) \Delta X_s\Big) \cdot 1_{\roi{0}{\zeta}}  \\
X(\kappa)        &:= (X-\check{X}(\kappa)).
\end{align}
The jumps $\Delta X(\kappa)=\kappa(\Delta X)\Delta X$ are bounded on $\roi{0}{\zeta}$. Hence by Lemma \ref{lem:boundedjumps} $X(\kappa)^{(j)}\in\cs_p^\dagger$ for $j=1,...,d$. Therefore there exists a unique decomposition
\begin{align} \label{specialdecomposition}
  X(\kappa)^{(j)}=X_0^{(j)}+M(\kappa)^{(j)}+B(\kappa)^{(j)}
\end{align}
where $M(\kappa)^{(j)}\in \cl^\dagger$ and $B(\kappa)^{(j)}$ is a predictable process in $\cv^\dagger$. The first characteristic of $X$ is defined to be $B:=B(\kappa)=(B(\kappa)^{(1)},...,B(\kappa)^{(d)})'$.

The localized second characteristic $C(n)$ is defined componentwise on $\roi{0}{\tau_n}$ for every $n\in\bbn$ by
\[
C^{jk}(n):=\big< X^{\tau_n,(j),c},X^{\tau_n,(k),c} \big>.
\]
For $m\geq n$ we obtain
\begin{align*}
[X^{\tau_m},X^{\tau_m}]_{t\wedge \tau_n}= [X,X]^{\tau_m}_{t\wedge \tau_n}=[X,X]^{\tau_m \wedge \tau_n}_t =[X,X]_t^{\tau_n} =[X^{\tau_n},X^{\tau_n}]_t
\end{align*}
and hence on $\roi{0}{\tau_n}$
\[
\Big<X^{\tau_m},X^{\tau_m}\Big>^c=\Big[X^{\tau_m},X^{\tau_m}\Big]^c =\Big[X^{\tau_n},X^{\tau_n}\Big]^c=
\Big<X^{\tau_n},X^{\tau_n}\Big>^c.
\]
Therefore we can set $C:=C(n)$ on $\roi{0}{\tau_n}$ for every $n\in\bbn$ and $C^{jk}:=\gamma$ on $\roi{\zeta}{\infty}$. $C$ is called the second characteristic of $X$.

The third characteristic $\nu$ is defined analogously as in the classical case, that is, $\nu$ is the compensator of the random measure of jumps $\mu^X$ restricted to $\roi{0}{\zeta}$:
\[
\mu^X(\omega;dt,dy):=\sum_{s < \zeta} 1_{\{\Delta X_s(\omega)\neq 0\} } \varepsilon_{(s,\Delta X_s(\omega))}(dt,dy)
\]
(cf. \cite{jacodshir} Section II.1) where $\varepsilon$ denotes the Dirac measure. On $\roi{\zeta}{\infty}$ $\nu$ is a priori not defined. For technical reasons we set for $H:\Omega\times [0,\infty[ \times (U-U)\to [0,\infty[$ which is $\cp\times \cb(U-U)$-measurable
\begin{align}
H(\omega,t,y)*\nu(\omega,t,y)= \gamma \text{ on } \roi{\zeta}{\infty}.
\end{align}

In Section 3 we use the following characterization of the characteristics (cf. \cite{jacodshir} introduction to Section II.2.d): Let
\begin{align} \label{charniceversion}
\begin{split}
B^{(j)}&=b^{(j)} \bullet F \\
C^{jk}&= c^{jk} \bullet F \\
\nu(\omega;dt,dy)&= K(\omega,t;dy) \, dF_t(\omega)
\end{split}
\end{align}
where \newline
(i) $F$ is a predictable process in $\ca_{loc}^{+,\dagger}$.\newline
(ii) $b=(b^{(j)})_{j\leq d}$ is a $d$-dimensional predictable process. \newline
(iii) $c=(c^{jk})_{j,k\leq d}$ is a predictable process with values in the set of all positive semidefinite $d\times d$- matrices. \newline
(iv) $K(\omega,t;dy)$ is a transition kernel from $(\Omega\times \bbr_+, \cp)$ into $(U-U,\cb(U-U))$ which satisfies on $\roi{0}{\zeta}$
\begin{align*}
  &K(\omega,t;\{0\})=0, \hspace*{10mm} \int_{U-U} \big( \abs{y}^2 \wedge 1 \big) K(\omega,t;dy) \leq 1 \\
  &\Delta F_t(\omega) > 0 \Rightarrow b_t(\omega) = \int_{U-U} y\kappa(y)K(\omega,t;dy) \\
  &\Delta F_t(\omega) K(\omega,t;{U-U})\leq 1.
\end{align*}
\emph{Remark: }If $(B,c,\nu)$ are the characteristics of a semimartingale with killing one can always find a version of the characteristics admitting a representation as above. To this end a localized reformulation of \cite{jacodshir} Proposition 2.9 is used. These `good versions' of the characteristics satisfy $\big(\abs{y}^2\wedge 1\big) * \nu(\cdot;[0,t]\times dy) \in \ca_{loc}^\dagger$ and identically on $\roi{0}{\zeta}$:
\begin{align}
&s\leq t \Rightarrow (C_t^{jk}-C_s^{jk})_{j,k\leq d} \text{ is a positive semidefinite matrix}; \\
&\nu(\cdot;\{t\}\times (U-U) ) \leq 1; \\
&\Delta B_t = \int_{U-U} y \kappa(y) \, \nu(\cdot; \{t\} \times dy).
\end{align}
We define
\[
F(\xi):=i B_t'\xi -\frac{1}{2} \xi'C_t \xi + \int_{U-U} \Big( e^{iy'\xi}-1-iy'\xi\kappa(y)\Big) \nu(\cdot;[0,t],dy)
\]
Then $F(\xi)_t = a(\xi)_t \bullet F_t$ where
\[
a(\xi)_t = i b_t'\xi - \frac{1}{2} \xi'c_t \xi + \int_{U-U} \Big( e^{iy'\xi} -1-iy'\xi \kappa(y) \Big) K(\cdot,t;dy).
\]
Therefore $(F(\xi)_t)_{t\geq 0}$ is a complex valued process with killing which is predictable and of finite variation.

\begin{Thm} \label{thm:charchar}
Let $X$ be an adapted c\`adl\`ag process killed at $\zeta$. Furthermore let $b,c,K,F$ be such that (i)-(iv) above are satisfied and let $(B,C,\nu)$ be given by \eqref{charniceversion}. The following are equivalent:
\begin{itemize}
  \item[(i)] $X\in\vec{\cs}^\dagger$ and it admits the characteristics $(B,C,\nu)$.
  \item[(ii)] For every $\xi\in\bbr^d$ the process given by
    \[
      e^{i\xi'X}-(e^{i\xi'X_-})\bullet F(\xi)
    \]
    is a complex valued local martingale with killing, that is, its real and imaginary part are in $\cm_{loc}^\dagger$.
  \item[(iii)] For every $f:\bbr^d_\Delta\to \bbr$ such that $f|_{U}\in C_b^2(U)$, $f|_{U_\Delta^c}=0$ and $f(\Delta)=\gamma$ the process
    \begin{align*}
      f(X_t)-f(X_0)&-\sum_{j\leq d} \partial_j f(X_t-) \bullet B_t^{(j)}
                  -\frac{1}{2} \sum_{j,k\leq d} \partial_j \partial_k f(X_t-) \bullet C_t^{jk}\\
                 &-\Big( f(X_{t-}+y)-f(X_{t-})-\sum_{j\leq d} \partial_j f(X_{t-}) (y^{(j)}) \kappa(y) \Big) * \nu(\cdot;dt,dy)
    \end{align*}
    is a local martingale with killing.
\end{itemize}
\end{Thm}

In order to prove this we need the following lemma which is taken from \cite{bergforst} (Theorem 10.8).

\begin{Lem} \label{lem:tripletunique}
Let $\ell\in\bbr^d$, $Q$ be a positive semidefinite $d\times d$-matrix and $N$ be a L\'evy measure on $U-U$. Then the representation of the function $\psi:\bbr^d \to \bbc$ given by
\[
  \psi(\xi)=i\ell'\xi -\frac{1}{2} \xi'Q\xi + \int_{U-U} \Big( e^{iy'\xi}-1-iy'\xi \kappa(y) \Big) \, N(dy) \hspace{10mm} (\xi\in\bbr^d)
\]
with such $(\ell, Q, N(dy))$ is unique.
\end{Lem}

\textbf{Proof of the theorem:}
`$(i)\Rightarrow (iii)$' First we use the decomposition \eqref{specialdecomposition} and write $M:=M(\kappa)$ and $B:=B(\kappa)$. For $n\in\bbn$ we use It\^o's formula and obtain
\begin{align*}
f(X^{\tau_n})-f(X_0)&=\sum_{j\leq d} \partial_j f(X_-^{\tau_n}) \bullet M^{\tau_n,(j)} + \sum_{j\leq d} \partial_j f(X_-^{\tau_n}) \bullet B^{\tau_n,(j)} \\
&\hspace{4mm}+\sum_{s\leq \cdot} \sum_{j\leq d} \partial_j f(X_{s-}^{\tau_n}) \Big(\Delta X_s^{\tau_n,(j)}-\kappa(\Delta X_s^{\tau_n})(\Delta X_s^{\tau_n})^{(j)}\Big) \\
&\hspace{4mm}+\frac{1}{2} \sum_{k,j\leq d} \partial_j \partial_k f(X_-^{\tau_n}) \bullet \Big< X^{\tau_n, (j), c}, X^{\tau_n, (k), c} \Big>\\
&\hspace{4mm}+\sum_{s\leq \cdot} \Bigg( f(X_s^{\tau_n})-f(X_{s-}^{\tau_n})-\sum_{j\leq d} \partial_j f(X_{s-}^{\tau_n}) \Delta X_s^{\tau_n,(j)} \Bigg) \\
&= \sum_{j\leq d} \partial_j f(X_-^{\tau_n}) \bullet M^{\tau_n,(j)}+ \sum_{j\leq d} \partial_j f(X_-^{\tau_n}) \bullet B^{\tau_n,(j)} \\
&\hspace{4mm}+ \frac{1}{2} \sum_{j,k\leq d} \partial_j \partial_k f(X_-^{\tau_n}) \bullet C^{\tau_n,jk} \\
&\hspace{4mm}+ H(t,y) 1_{\ci{0}{\tau_n}}* \mu^X(\cdot;dt,dy)
\end{align*}
where $H(t,y):=f(X_{t-}^{\tau_n}+y)-f(X_{t-}^{\tau_n}) - \sum_{j\leq d} \partial_j f(X_{t-}^{\tau_n}) y^{(j)}\kappa(y)$.  Literally as in the classical case (cf. \cite{jacodshir} Theorem II.2.42) one obtains that the last term on the right-hand side is in $\ca_{loc}$. The result follows by \cite{jacodshir} Proposition II.1.28 and the fact that $(W  1_{\ci{0}{\tau_n}})*(\mu^X-\nu) = W\bullet (1_{\ci{0}{\tau_n}}*(\mu^X-\nu))$.

`$(iii)\Rightarrow (ii)$' Use (iii) on the function
\[
f(y)= \begin{cases} e^{iy'\xi} &\text{if } y\in U \\ \gamma &\text{if } y=\Delta \\
                    0 & \text{else}\end{cases}
\]
respectively on its real and imaginary part. The claim is then proved as in the classical case.

`$(ii)\Rightarrow (i)$' Using the hypothesis and the fact that $F(\xi)$ is of finite variation, we obtain that $e^{i\xi'X}$ is a complex valued semimartingale with killing for every $\xi\in\bbr^d$. Therefore $\sin(rX^{(j)})\in\cs^\dagger$ for each $r\in\bbr$. There exists a function $f\in C^2(\bbr)$ such that $f(\sin(x))=x$ on $\abs{x}\leq 1/2$. We set $\sigma_n:=\inf \{t\geq 0:  X_t= \Delta \text{ or } |X_t^{(j)}|>n/2 \}$ and $\rho_n:=\sigma_n \wedge \tau_n$. Then we have
\[
X^{(j)}=nf\Big(\sin\big(X^{(j)}/n \big)\Big) \hspace{1cm}\text{ on } \roi{0}{\rho_n}
\]
and $\rho_n\uparrows \zeta$. By Theorem \ref{thm:prelocal} we obtain $X^{(j)}\in\cs^\dagger$. \newline
Now let $(\widetilde{B},\widetilde{C},\widetilde{\nu})$ be good versions of the characteristics of $X$. (see the remark above). Associate the following process with these characteristics
\[
\widetilde{F}(\xi)=i\xi \widetilde{B}_t -\frac{1}{2} \xi'\widetilde{C}_t \xi + \int \Big( e^{iy'\xi}-1-iy'\xi\kappa(y)\Big) \widetilde{\nu}(\cdot;[0,t],dy)
\]
This process is again killed at $\zeta$. We already know that $(i)$ implies $(ii)$. Therefore
\begin{align} \label{Frepunique}
e^{i\xi'X}-\Big(e^{i\xi'X_-} \bullet \widetilde{F}(\xi)\Big) \in \cl^\dagger
\end{align}
where the second term is a predictable process in $\cv^\dagger$. Since $e^{i\xi'X}$ is bounded on $\roi{0}{\zeta}$ it is a special semimartingale with killing. Hence the decomposition \eqref{Frepunique} is unique and therefore
\[
e^{i\xi'X_-} \bullet F(\xi) = e^{i\xi'X_-} \bullet \widetilde{F}(\xi).
\]
Integrating $e^{-i\xi'X_-}$ against both sides we obtain $F(\xi)=\widetilde{F}(\xi)$, that is, the set $N$ of all $\omega$ for which there exist $\xi\in\bbq^d$ and $t\in\bbq_+$ such that $F(\xi)_t(\omega)\neq \widetilde{F}(\xi)_t(\omega)$ is $\bbp$-Null. On $\roi{\zeta}{\infty}$ all the processes $X,B,C,F(\xi), \widetilde{B}, \widetilde{C}, \widetilde{F}(\xi)$ as well as $\nu(\cdot,[0,t],E),\widetilde{\nu}(\cdot,[0,t],E)$ for $E\in \cb(U-U)$ are equal to $\Delta$. On $\roi{0}{\zeta}$ we argue as follows: the function $\psi$ in Lemma \ref{lem:tripletunique} is continuous and hence fully characterized by its values on $\bbq^d$. On the complement of $N$ we have
\[
B_t=\widetilde{B}_t, C_t=\widetilde{C}_t, \nu(\cdot;[0,t],E)=\widetilde{\nu}(\cdot;[0,t],E) \hspace{1cm}\text{ for all } t\in\bbq_+, E\in\cb(U-U).
\]
Since the processes are right continuous the equalities hold on $\bbr_+$. Hence the result.
\hfill $\square$

\begin{Def} \label{def:ito}
Let $\textbf{X}=(\Omega, \cf ,(\cf_t)_{t\geq 0},(X_t)_{t\geq 0},\bbp^x)_{x\in U_\Delta}$ be a Markov process on $(U,\cb(U))$ with killing. We call \textbf{X} \emph{It\^o process with killing} if $\textbf{X}\in\vec{\cs}^\dagger$ and its semimartingale characteristics $(B,C,\nu)$ with respect to the fixed cut-off function $\kappa$ can be written as
\begin{align} \begin{split}
B_t^{(j)}(\omega)&=\int_0^t  \ell^{(j)}(X_s(\omega)) \, dF_s \\
C_t^{jk}(\omega)&= \int_0^t Q^{jk}(X_s(\omega)) \, dF_s \\
\nu(\omega;ds,dy)&=N(X_s(\omega),dy) \, dF_s
\end{split}
\end{align} for every $\bbp^x, (x\in U)$ where for every $z\in U$: $\ell(z) \in \bbr^d$, $Q(z)$ is a positive semidefinite matrix, $N(z,dy)$ is a Borel transition kernel such that $N(z,\{0\})=0$ and $F_s=s\cdot 1_{\roi{0}{\zeta}}(s)+ \Delta \cdot 1_{\roi{\zeta}{\infty}}(s)$.
\end{Def}

\emph{Remark:} It\^o processes in the above sense (on $\bbr^d$ without killing) were introduced in \cite{vierleute} as a sub-class of Hunt semimartingales. In \cite{cinlarjacod81} it is shown that they can be characterized as a class of solutions of certain SDEs. The stochastic symbol for this class of processes was introduded in \cite{mydiss}.

\section*{Acknowledgements}

This paper contains some of the results of my Ph.D thesis written under the guidance of Professor Ren\'e L. Schilling, to whom I am deeply grateful for many helpful suggestions. Furthermore I would like to thank an anonymous referee who took a lot of time in very careful proof reading and in writing an extraordinarily detailed report. His/her comments helped to improve the paper. Financial support of the DFG (German Science Foundation) for the project SCHN 1231/1-1 is gratefully acknowledged. 


\bibliographystyle{plainlf}

\end{document}